\documentclass[11pt,a4paper]{article}
\usepackage[T1]{fontenc}
\usepackage{gfsdidot}
\usepackage{authblk}
\usepackage{marvosym}
\usepackage{graphicx}
\usepackage{enumitem}
\usepackage{apalike}
\usepackage{vmargin}
\setmarginsrb{2cm}{2cm}{2cm}{2cm}{0mm}{0mm}{0mm}{10mm}
\usepackage[colorlinks=true, allcolors=blue, breaklinks=true]{hyperref}

\title{\vspace*{-1.25cm} \bfseries Active participation and student journal in\\ Confucian heritage culture mathematics classrooms}
\author{\normalsize Natanael Karjanto\thanks{\Letter: \url{natanael@skku.edu} \href{https://orcid.org/0000-0002-6859-447X}{\includegraphics[scale=0.08]{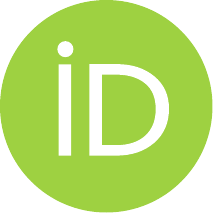}}}}
\affil{Department of Mathematics, University College, Natural Science Campus\\ Sungkyunkwan University, Suwon~16419, Republic of Korea}

\date{\vspace*{-0.5cm} \scriptsize Updated \today}

\begin{document}
\maketitle
\thispagestyle{empty}
\vspace*{-0.6cm}

\begin{abstract}
\noindent
This article discusses an effort to encourage student-instructor interactive engagement through active learning activities during class time. We not only encouraged our students to speak out when an opportunity arises but also required them to record their active participation in a student journal throughout the semester. In principle, any activities which constitute active learning can and should be recorded in a \emph{student journal}. These include, but are not limited to, reading definition, theorem, problem, etc.; responding to questions and inquiries; asking questions; and pointing out some mistakes during class time. Despite an incentive for this participation, our experience teaching different mathematics courses in several consecutive semesters indicates that many students resist speaking out publicly, submitting empty journals at the end of the semester instead. Students' feedback on teaching evaluation at the end of the semester reveals that many students dislike and are against the idea of active participation and recording it in the journal. This paper discusses the reason behind this resistance and provides some potential remedies to alleviate the situation. \\

\end{abstract}

\section{INTRODUCTION} 

\subsection{Literature study}
Adopting and implementing interactive engagement and students' active participation are not only fostering a fun and interesting environment for studying but also cultivate students with the goal of deeper learning and long-term retention rate. The literature contains a non-exhaustive list of pedagogy involving some kind of active learning in various disciplines. There exists ample evidence that active learning is superior to a passive one, from improving students' grades, minimizing the number of failures and increasing the lifespan of retention learning rate~\cite{Crouch01,Freeman14}.

We observed two potential challenges which may limit the success of active learning strategy: Confucian Heritage Culture (CHC) and English-Medium Instruction (EMI). The former is not only known for its quiet class culture and passive-receptive learning style but also includes cultural differences and language barriers. Hence, both challenges are closely related. Indeed, typical CHC students are often reluctant in expressing opinions and are not accustomed to classroom participation. Some of them even go further by resisting curriculum, pedagogy, and context, as recently revealed by~\cite{Huang18}.

EMI is not only popular in Europe but also in Asia because offering it may improve the university's reputation and ranking~\cite{Cho12}. This ambition can be hampered due to the lack of language proficiency, from both instructors and students alike, as evidenced in the Korean context~\cite{Byun11}. We should learn from our European colleagues how to implement EMI successfully~\cite{Klaassen01}. Admittedly, the issue surrounding EMI is not simple and any bold step in adopting it has fostered debates, challenges, and controversies~\cite{Doiz12}.

\subsection{Theoretical framework}
The theoretical framework for this study is a combination of journal writing, CHC classroom, and EMI environment. For the student journal, it is a tool to acquire and improve reflective thinking~\cite{Dewey33,Schon87}. For CHC, it is based on both the conventional teacher-centered approach to teaching and learning as well as the `activity theory' which enhances constructivist learning reforms~\cite{Engestrom99,Ng09}. For EMI, it is based on the `dynamic bilingual education', which involves the practice of using English for communication purposes as well as promoting multicultural awareness~\cite{Baker11,Garcia09}.

\subsection{Research question}
We are interested in answering the following research questions:
\begin{enumerate}[leftmargin=1.3em]
\item What is students' feedback after experiencing active learning for which they also need to record their participation in a student journal?
\item What are some challenges and potential remedies for this type of pedagogy, particularly in the context of the CHC-EMI environment?
\end{enumerate}

\section{METHODOLOGY}
We only consider the qualitative aspect of active participation and student journal using the students' feedback and perception collected at the end of the semester.

\subsection{Participant}
The participants in this study are the students who are enrolled in three mathematics courses offered at the Natural Science Campus of our university: Single Variable Calculus (two sections, 134 students), Multivariable Calculus (two sections, 84 students) and Linear Algebra (six sections, 336 students), for two years from Fall 2016 until Spring 2018. The total number of participants is 554 and their age ranges from 18 to 24 years old. We also adopted the convenience sampling method due to its efficiency and accessibility.

\subsection{Measurement}
We obtain the students' feedback from the online questionnaire administered by the Academic Affairs Team. We focus on the second part of this questionnaire which solicited the students' suggestions for the instructor in improving the teaching.

\section{RESULT}
The submission rate of the journal ranges from 38\% to 94\%, depending on the course and the semester. From those who submitted the student journal, around 30\% of the students have never participated during the class, as indicated by empty journal submission or unrelated comments were written down just for the sake of filling out the form. From our class observation, only around 10\% of the students participated regularly and filled out the student journal. Digging in further on students' feedback confirms our initial hypothesis that many students dislike the idea of active participation and recording it in a student journal, as the following comments show.\\
{\small ``I think the student journal doesn't improve anything. Hope we don't do it.'' (Spring 2017, Single-Variable Calculus)} \\
{\small ``I hate commenting.'' (Fall 2016, Multivariable Calculus)} \\
{\small ``I think Student journal can't help my LA studying.'' (Fall 2016, Linear Algebra)} \\
{\small ``I think student journals are not useful.'' (Fall 2017, Linear Algebra)}   \\
{\small ``It seems that the instructor was enforcing too much on active participation.'' (Spring 2018, Linear Algebra)}

From these comments, we observe that generally the students do not like the idea of participating in the class, particularly when they have to speak out publicly in front of their peers and write down what they have done in the student journal. Many students prefer the passive-receptive learning style by observing and listening to the instructor talking in front of the classroom. 

\section{DISCUSSION}

Several factors may contribute to the students' resistance in participating actively in class. During their previous educational experience, many of them grow up in a learning environment where active participation is not required and the learning style is passive-receptive. In several cases, speaking up publicly is frowned upon since the culture dictates them to remain quiet from an early age, to listen attentively in class, particularly when the teacher is talking and explaining, and not to challenge authority even though there exist some obvious mistakes. This type of learning style is common among students who were educated in the CHC environment. In each student's mind, there is a `little-Confucius' subconsciously reminding them to do what they are used to doing, despite an invitation to do otherwise. 

In addition to the cultural background, language skill also plays an important role. The mathematics courses are offered as EMI and the students have diverse English proficiency, as well as mathematical ability. From our observation, many students in the West generally have higher English proficiency even though English is not their first language in comparison to the students. Many factors can contribute to this, including how similar the students' first language is to English, the quality of English education at the secondary level, and admission criteria of English proficiency, such as TOEFL or IELTS score. The majority of students are selected through the national examination. English testing is a component, but communication skills and particularly public speaking are not examined.

Regarding the diversity of mathematical ability, our observation indicates that academically-prepared students tend to cooperate more in active participation and in writing the student journal, even though some of them have limited English ability. It also makes sense that weaker students tend to avoid participation because they might not understand the material or get lost in the discussion.

The remedy for the latter is beyond our discussion in this paper, but for the former, we could address the culture before implementing active learning pedagogy. The purpose is not entirely to override the culture as a whole but to have a common ground of learning culture where we stand at the same level and could progress together to make successful teaching and learning. We may explain to the students at the beginning of every semester that active participation is required if they wish to enroll in a particular course and that there will be points awarded for active participation, and these will be counted toward their final score. Furthermore, we can assure them that it is all right to be noisy during the class because it is part of the learning process. We conjecture that the pedagogy of interactive engagement in this special CHC-EMI environment will be more effective if both the instructors and the students possess high proficiency not only in English but also in Korean.

\subsection*{Acknowledgment}

The author gratefully acknowledges Dr. Kristiana Wijaya for the invitation to present as one of the invited speakers at the International Conference on Mathematics, Geometry, Statistics, and Computation (IC-MaGeStiC 2021), held on 27 November 2021 at Jember University, East Java, Indonesia.

\vfill
\noindent
{\small How to cite this article:\\
Karjanto, N. (2022). Active participation and student journal in Confucian heritage culture mathematics classrooms. In Wijaya, K. (Editor) \textit{Proceedings of the International Conference on Mathematics, Geometry, Statistics, and Computation (IC-MaGeStiC 2021), Advances in Computer Science Research}, volume 96, pages~89--91. Dordrecht, the Netherlands: Atlantis Press. \hfill \url{https://arxiv.org/abs/1912.07837} \\ \url{https://doi.org/10.2991/acsr.k.220202.018}. \par}
\end{document}